\documentclass[12pt]{article}
\usepackage{amsmath,amssymb,amsthm}

\def\barr{\begin{array}}
\def\earr{\end{array}}

\title{The number of cyclic subgroups of finite abelian groups and Menon's identity}
\author{Marius T\u arn\u auceanu}
\date{October 1, 2018}

\begin{document}

\maketitle

\begin{abstract}
In this note, we give a new formula for the number of cyclic subgroups of a finite abelian group. This is based on applying
the Burnside's lemma to a certain group action. Also, it generalizes the well-known Menon's identity.
\end{abstract}

{\small
\noindent
{\bf MSC 2010\,:} Primary 11A25; Secondary 20D60, 20K01.

\noindent
{\bf Key words\,:} Menon's identity, Burnside's lemma, group action, cyclic subgroup, abelian group, power automorphism group.}

\section{Introduction}

One of the most interesting arithmetical identities is due to P.K. Menon \cite{9}.

\smallskip\noindent{\bf Menon's identity.}
{\it For every positive integer $n$ we have
\begin{equation}
\sum_{a\in\mathbb{Z}^*_n}\gcd(a-1,n)=\varphi(n)\tau(n),\nonumber
\end{equation}where $\mathbb{Z}^*_n$ is the group of units of the ring $\mathbb{Z}_n=\mathbb{Z}/n\mathbb{Z}$, $\gcd(,)$ represents
the greatest common divisor, $\varphi$ is the Euler's totient function and $\tau(n)$ is the number of divisors of $n$.}
\smallskip

This identity has many generalizations derived by several authors (see, for example, \cite{1}-\cite{8} and \cite{10,11,14,15,17,19,20}). A usual technique to prove results of this type is based on the Burnside's lemma (also called the Cauchy-Frobenius lemma - see \cite{12}) concerning group actions.

\smallskip\noindent{\bf Burnside's lemma.}
{\it Let $G$ be a finite group acting on a finite set $X$ and
\begin{equation}
Fix(g)=\{x\in X \mid g\circ x=x\}, \mbox{ for all } g\in G.\nonumber
\end{equation}Then the number of distinct orbits is
\begin{equation}
N=\frac{1}{|G|}\sum_{g\in G}|Fix(g)|.
\end{equation}}
\smallskip

In the following let $G$ be a finite abelian group of order $n$ and $$G=G_1\times\cdots\times G_k$$be the primary decomposition of $G$, where $G_i$ is a $p_i$-group for all $i=1,\dots,k$. Then every $G_i$ is of type $$G_i=\mathbb{Z}_{p_i^{\alpha_{i1}}}\times\cdots\times\mathbb{Z}_{p_i^{\alpha_{ir_i}}},$$where $1\leq\alpha_{i1}\leq\cdots\leq\alpha_{ir_i}$. We will apply the Burnside's lemma to the natural action of the power automorphism group ${\rm Pot}(G)$ on $G$. Recall that an automorphism $f$ of $G$ is called a \textit{power automorphism} if $f(H)=H, \,\forall\, H\leq G$, and that the set ${\rm Pot}(G)$ of all power automorphisms of $G$ is a subgroup of ${\rm Aut}(G)$. Also, it is well-known that every power automorphism of a finite abelian group is \textit{universal}, i.e. there exists an integer $m$ such that $f(x)=mx$ for all $x\in G$. The structure of ${\rm Pot}(G)$ is given by Theorem 1.5.6 in \cite{13}:
\begin{equation}
{\rm Pot}(G)\cong{\rm Pot}(G_1)\times\cdots\times{\rm Pot}(G_k)\cong{\rm Aut}(\mathbb{Z}_{p_1^{\alpha_{1r_1}}})\times\cdots\times{\rm Aut}(\mathbb{Z}_{p_k^{\alpha_{kr_k}}}).
\end{equation}

Our main result is stated as follows.

\smallskip\noindent{\bf Theorem 1.}
{\it Under the above notations, we have
\begin{equation}
\prod_{i=1}^k\!\!\!\sum_{^{\ \ \ 1\leq m_i\leq p_i^{\alpha_{ir_i}}}_{\ \ \ \ p_i\nmid m_i}}\prod_{j=1}^{r_i}\gcd(m_i-1,p_i^{\alpha_{ij}})=\varphi(\exp(G))|L_1(G)|,
\end{equation}where $\exp(G)$ is the exponent of $G$ and $|L_1(G)|$ is the number of cyclic subgroups of $G$.}
\smallskip

Clearly, (3) gives a new formula to compute the number of cyclic subgroups of a finite abelian group (for other such formulas see e.g. \cite{16,18}). We exemplify it in a particular case.

\smallskip\noindent{\bf Example.} The finite abelian group $$G=\mathbb{Z}_2\times\mathbb{Z}_{12}\times\mathbb{Z}_{72}\cong\left(\mathbb{Z}_2\times\mathbb{Z}_{2^2}\times\mathbb{Z}_{2^3}\right)\times\left(\mathbb{Z}_3\times\mathbb{Z}_{3^2}\right)$$has
$\exp(G)=2^33^2=72$ and so $\varphi(\exp(G))=\varphi(72)=24$. Then (3) leads to
$$|L_1(G)|=\frac{1}{24}\left(\sum_{^{\!\!\ \ 1\leq m_1\leq 2^3}_{\ \ \ 2\nmid m_1}}\prod_{j=1}^3\gcd(m_1-1,2^{\alpha_{1j}})\right)\left(\sum_{^{\!\!\ \ 1\leq m_2\leq 3^2}_{\ \ \ 3\nmid m_2}}\prod_{j=1}^2\gcd(m_2-1,3^{\alpha_{2j}})\right)\vspace{-5mm}$$
\begin{align*}
&\hspace{13,5mm}=\frac{1}{24}\left(\gcd(0,2^1)\gcd(0,2^2)\gcd(0,2^3)+\gcd(2,2^1)\gcd(2,2^2)\gcd(2,2^3)+\right.\\
&\hspace{18mm}\left.\gcd(4,2^1)\gcd(4,2^2)\gcd(4,2^3)+\gcd(6,2^1)\gcd(6,2^2)\gcd(6,2^3)\right)\\
&\hspace{18mm}\left(\gcd(0,3^1)\gcd(0,3^2)+\gcd(1,3^1)\gcd(1,3^2)+\gcd(3,3^1)\gcd(3,3^2)+\right.\\
&\hspace{18mm}\left.\gcd(4,3^1)\gcd(4,3^2)+\gcd(6,3^1)\gcd(6,3^2)+\gcd(7,3^1)\gcd(7,3^2)\right)\\
&\hspace{13,5mm}=\frac{1}{24}\left(64+8+32+8\right)\left(27+1+9+1+9+1\right)\\
&\hspace{13,5mm}=224.
\end{align*}

We remark that if the group $G$ is cyclic of order $n$, then $r_i=1,\, \forall\,i=1,...,k$, $\exp(G)=p_1^{\alpha_{11}}\cdots p_k^{\alpha_{k1}}=n$ and $|L_1(G)|=\tau(n)$. Thus the equality (3) becomes
\begin{equation}
\prod_{i=1}^k\!\!\!\sum_{^{\ \ \ 1\leq m_i\leq p_i^{\alpha_{i1}}}_{\ \ \ \ p_i\nmid m_i}}\gcd(m_i-1,p_i^{\alpha_{i1}})=\varphi(n)\tau(n).
\end{equation}Since
\begin{equation}
\mathbb{Z}^*_{p_1^{\alpha_{11}}}\times\cdots\times\mathbb{Z}^*_{p_k^{\alpha_{k1}}}\cong\mathbb{Z}^*_n,\nonumber
\end{equation}(4) can be rewritten as
\begin{equation}
\sum_{^{ \ \ \ 1\leq m\leq n}_{ \ \ \gcd(m,n)=1}}\gcd(m-1,n)=\varphi(n)\tau(n),\nonumber
\end{equation}that is we obtained the Menon's identity.
\smallskip

Two immediate consequences of Theorem 1 are the following.

\smallskip\noindent{\bf Corollary 2.}
{\it Let $m$ and $n$ be two positive integers, $l={\rm lcm}(m,n)$, and $p_1,...,p_k$ be the primes dividing $l$. Write $m=p_1^{\alpha_1}\cdots p_k^{\alpha_k}$ and $n=p_1^{\beta_1}\cdots p_k^{\beta_k}$, where $\alpha_i$ and $\beta_i$ can be possibly zero. Then
\begin{equation}
|L_1(\mathbb{Z}_m\times\mathbb{Z}_n)|=\frac{1}{\varphi(l)}\,\prod_{i=1}^k\!\!\!\!\!\!\!\!\!\!\!\!\!\!\!\sum_{^{\ \ \ \ \ \ \ 1\leq m_i\leq p_i^{\max\{\alpha_i,\beta_i\}}}_{\ \ \ \ \ \ \ \ p_i\nmid m_i}}\!\!\!\!\!\!\gcd(m_i-1,p_i^{\alpha_i})\gcd(m_i-1,p_i^{\beta_i}).
\end{equation}}

\smallskip\noindent{\bf Corollary 3.}
{\it Let $n$ be a positive integer and $n=p_1^{\alpha_1}\cdots p_k^{\alpha_k}$ be the decomposition of $n$ as a product of prime factors. Then for every $r\in\mathbb{N}^*$ we have
\begin{equation}
|L_1(\mathbb{Z}_n^r)|=\frac{1}{\varphi(n)}\,\prod_{i=1}^k\!\!\!\!\sum_{^{\ \ \ 1\leq m_i\leq p_i^{\alpha_i}}_{\ \ \ \ p_i\nmid m_i}}\!\!\!\!\!\!\gcd(m_i-1,p_i^{\alpha_i})^r.
\end{equation}}
\smallskip

Note that (4) can be obtained from (5) or (6) by taking $m=1$ or $r=1$, respectively. Thus, these equalities can be also seen as generalizations of Menon's identity.

\section{Proof of Theorem 1}

The natural action of ${\rm Pot}(G)$ on $G$ is
$$f\circ a=f(a),\, \forall\, (f,a)\in{\rm Pot}(G)\times G.$$By using the direct decompositions of ${\rm Pot}(G)$ and $G$ in Section 1, it can be written as
$$(f_1,...,f_k)\circ (a_1,...,a_k){=}(f_1(a_1),...,f_k(a_k)),\, \forall\, (f_i,a_i)\in{\rm Pot}(G_i)\times G_i, i=1,...,k.$$

First of all, we will prove that two elements $a,b\in G$ are contained in the same orbit if and only if they generate the same cyclic subgroup of $G$. Indeed, if $a$ and $b$ belong to the same orbit, then there exists $f\in{\rm Pot}(G)$ such that $b=f(a)$. Since $f$ is universal, it follows that $b=ma$ for some integer $m$. Then $b\in\langle a\rangle$, and so $\langle b\rangle\subseteq\langle a\rangle$. On the other hand, since a group automorphism preserves the element orders, we have $o(a)=o(b)$. Therefore $\langle a\rangle=\langle b\rangle$. Conversely, assume that $\langle a\rangle=\langle b\rangle$, where $a=(a_1,...,a_k)$ and $b=(b_1,...,b_k)$. Then $\langle a_i\rangle=\langle b_i\rangle$, $\,\forall\, i=1,...,k$. This implies that for every $i$ there is an integer $m_i$ such that $b_i=m_ia_i$ and $\gcd(m_i,o(a_i))=1$. Remark that if $a_i=1$ then we must have $m_i=1$, while if $a_i\neq 1$ then $\gcd(m_i,p_i)=1$. Consequently, in both cases $p_i\nmid m_i$. This shows that the map
$$f_i:G_i\longrightarrow G_i, f_i(x_i)=m_ix_i,\, \forall\, x_i\in G_i$$is a power automorphism of $G_i$. Then $f=(f_1,...,f_k)\in{\rm Pot}(G)$ and $f(a)=b$, i.e. $a$ and $b$  are contained in the same orbit. Thus, the number of distinct orbits is $$N=|L_1(G)|.$$

Next we will focus on the right side of (1). Note that the group isomorphism (2) leads to
$$|{\rm Pot}(G)|=\prod_{i=1}^k|{\rm Aut}(\mathbb{Z}_{p_i^{\alpha_{ir_i}}})|=\prod_{i=1}^k\varphi(p_i^{\alpha_{ir_i}})=\varphi(\prod_{i=1}^k p_i^{\alpha_{ir_i}})=\varphi(\exp(G)).$$Also, for every $f=(f_1,...,f_k)\in{\rm Pot}(G)$ and every $a=(a_1,...,a_k)\in G$ we have
$$a\in Fix(f) \Leftrightarrow a_i\in Fix(f_i),\, \forall\, i=1,...,k,$$implying that $$|Fix(f)|=\prod_{i=1}^k |Fix(f_i)|.$$On the other hand, since
${\rm Pot}(G_i)\cong{\rm Aut}(\mathbb{Z}_{p_i^{\alpha_{ir_i}}})$, every $f_i$ is of type
$$f_i(x_i)=m_ix_i \mbox{ with } p_i\nmid m_i.$$Then for $x_i=(x_{i1},...,x_{ir_i})\in G_i$ we have
$$\hspace{-30mm}x_i\in Fix(f_i)\Leftrightarrow (m_i-1)x_{ij}=0 \mbox{ in } \mathbb{Z}_{p_i^{\alpha_{ij}}},\, \forall\, j=1,...,r_i$$
\begin{align*}
&\hspace{27mm}\Leftrightarrow p_i^{\alpha_{ij}}\mid (m_i-1)x_{ij},\, \forall\, j=1,...,r_i\\
&\hspace{27mm}\Leftrightarrow\frac{p_i^{\alpha_{ij}}}{\gcd(m_i-1,p_i^{\alpha_{ij}})}\mid x_{ij},\, \forall\, j=1,...,r_i\\
&\hspace{27mm}\Leftrightarrow x_{ij}=c\frac{p_i^{\alpha_{ij}}}{\gcd(m_i-1,p_i^{\alpha_{ij}})} \mbox{ with } c=0,...,\gcd(m_i-1,p_i^{\alpha_{ij}})-1,\\ &\hspace{34mm}\forall\, j=1,...,r_i,
\end{align*}and so $$|Fix(f_i)|=\prod_{j=1}^{r_i}\gcd(m_i-1,p_i^{\alpha_{ij}}).$$Thus, the right side of (1) becomes
$$\frac{1}{\varphi(\exp(G))}\sum_{f\in{\rm Pot}(G)}\!\!\!|Fix(f)|=\frac{1}{\varphi(\exp(G))}\sum_{f_1\in{\rm Pot}(G_1)}\!\!...\!\!\!\sum_{f_k\in{\rm Pot}(G_k)}|Fix(f_1)|\cdots|Fix(f_k)|$$
$$=\frac{1}{\varphi(\exp(G))}\prod_{i=1}^k\left(\sum_{f_i\in{\rm Pot}(G_i)}\!\!\!|Fix(f_i)|\right)=\frac{1}{\varphi(\exp(G))}\prod_{i=1}^k\!\!\!\!\!\sum_{^{\ \ \ 1\leq m_i\leq p_i^{\alpha_{ir_i}}}_{\ \ \ \ p_i\nmid m_i}}\prod_{j=1}^{r_i}\gcd(m_i-1,p_i^{\alpha_{ij}}),$$as desired.\hfill\rule{1,5mm}{1,5mm}

\vspace*{3ex}\small

\hfill
\begin{minipage}[t]{5cm}
Marius T\u arn\u auceanu \\
Faculty of  Mathematics \\
``Al.I. Cuza'' University \\
Ia\c si, Romania \\
e-mail: {\tt tarnauc@uaic.ro}
\end{minipage}

\end{document}